\documentclass[10pt]{article}
\usepackage{graphicx}
\newtheorem{theorem}{Theorem}
\newtheorem{proposition}{Proposition}
\newtheorem{lemma}{Lemma}
\newtheorem{corollary}{Corollary}

\input tcilatex

\begin{document}

\begin{titlepage}

\vskip0.5truecm

\vskip1.0truecm

\begin{center}

{\LARGE \bf Dynamics of homeomorphisms of the torus homotopic to Dehn twists}

\end{center}

\vskip  0.4truecm

\centerline {{\large Salvador Addas-Zanata, F\'abio A. Tal and Br\'aulio A. Garcia}}

\vskip 0.2truecm

\centerline { {\sl Instituto de Matem\'atica e Estat\'\i stica }}
\centerline {{\sl Universidade de S\~ao Paulo}}
\centerline {{\sl Rua do Mat\~ao 1010, Cidade Universit\'aria,}} 
\centerline {{\sl 05508-090 S\~ao Paulo, SP, Brazil}}
 
\vskip 0.7truecm

\begin{abstract}

In this paper we consider torus homeomorphisms $f$ homotopic to Dehn twists.
We prove that if the vertical rotation set of $f$ is reduced to zero, then  there exists
a compact connected essential "horizontal" set K, invariant under $f$.
In other words, if we consider the lift $\widehat{f}$ of $f$ to the cylinder, which
has zero vertical rotation number, then all points 
have uniformly bounded motion under iterates of $\widehat{f}$. Also, we give
a simple explicit condition which, when satisfied, implies that the vertical rotation set 
contains an interval and thus also implies positive topological entropy. 
As a corollary of the above results, we prove a version of Boyland's conjecture
 to this setting:
If $f$ is area preserving and has a lift $\widehat{f}$ to the cylinder 
with zero Lebesgue measure vertical rotation number, then either the orbits of all points are uniformly bounded under 
$\widehat{f}$, or there
are points in the cylinder with positive vertical velocity and others with negative
vertical velocity.

\end{abstract} 

\vskip 0.3truecm

\vskip 2.0truecm


\vskip 0.8truecm

\noindent{\bf e-mail:} sazanata@ime.usp.br, fabiotal@ime.usp.br and braulio@ime.usp.br

\vskip 1.0truecm


\vfill
\hrule
\noindent{\footnotesize{The first two authors are partially supported 
by CNPq, grants: 304803/06-5 and 304360/05-8. The third is supported by a
FAPESP grant 2008/10363-5.}}

\end{titlepage}

\baselineskip=6.2mm

\section{Introduction and main results}

In this paper we study homeomorphisms $f$ of the torus homotopic to Dehn
twists. These homotopy classes are in some way simpler to analyze than the
identity case. One of the reasons for this is the fact that there is no
sense in defining a two dimensional rotation set for torus maps homotopic to
Dehn twists, instead a vertical rotation set is defined, see expression (\ref
{rotset}).

Many important conjectures for homotopic to the identity maps have their
analogs in this setting. For instance, how is the rotation interval of a
minimal Dehn twist homeomorphism? Does the set of minimal Dehn twist $C^r$%
-diffeomorphisms $(r\geq 2)$ have no interior? If $f$ is a Dehn twist
homeomorphism which preserves area and has zero Lebesgue measure vertical
rotation number, is it true that either $f$ is more or less like an annulus
homeomorphism or the vertical rotation interval has no empty interior?

One of the main motivations for our work is a recent example of F. Tal and
A. Koropecki, where they present an area preserving torus homeomorphism $h$
homotopic to the identity, such that its rotation set is only $(0,0),$
satisfying the following property:

\begin{itemize}
\item  $h$ has a lift to the plane, denoted $\widetilde{h},$ such that $
\widetilde{h}$ has fixed points and some points in the plane have unbounded $
\widetilde{h}$-orbits in every direction
\end{itemize}

In other words, this example implies that the existence of sub-linear
displacement does not imply linear displacement, at least in the homotopic
to the identity class. In this work we show that maps homotopic to Dehn
twists have a different behavior. Before presenting our results, we need
some definitions.

\vskip0.2truecm

{\large {\bf Definitions:}}

\begin{enumerate}
\item  Let ${\rm T^2}={\rm I}\negthinspace {\rm R^2}/{\rm Z\negthinspace
\negthinspace Z^2}$ be the flat torus and let $p:{\rm I}\negthinspace {\rm %
R^2}\longrightarrow {\rm T^2}$ and $\pi :{\rm I}\negthinspace {\rm R^2}%
\longrightarrow S^1\times {\rm I}\negthinspace {\rm R}$ be the associated
covering maps. Coordinates are denoted as $(\widetilde{x},\widetilde{y})\in 
{\rm I}\negthinspace {\rm R^2,}$ $(\widehat{x},\widehat{y})\in S^1\times 
{\rm I}\negthinspace {\rm R}$ and $(x,y)\in {\rm T^2.}$

\item  Let $DT({\rm T^2})$ be the set of homeomorphisms of the torus
homotopic to a Dehn twist $(x,y)\longrightarrow (x+ky$ mod $1,y$ mod $1)$,
for some $k\in {\rm Z\negthinspace
\negthinspace Z^{*}}$, and let $DT(S^1\times {\rm I}\negthinspace {\rm R})$
and $DT({\rm I\negthinspace R^2})$ be the sets of lifts of elements from $DT(%
{\rm T^2})$ to the cylinder and plane. Homeomorphisms from $DT({\rm T^2})$
are denoted $f$ and their lifts to the vertical cylinder and plane are
respectively denoted $\widehat{f}$ and $\widetilde{f}.$

\item  Let $p_{1,2}:{\rm I}\negthinspace {\rm R^2}\longrightarrow {\rm I}%
\negthinspace {\rm R}$ be the standard projections; $p_1(\tilde x,\tilde
y)=\tilde x$ and $p_2(\tilde x,\tilde y)=\tilde y$. Projections on the
cylinder are also denoted by $p_1$ and $p_2.$

\item  Given $f\in DT({\rm T^2})$ and a lift $\widehat{f}\in DT(S^1\times 
{\rm I}\negthinspace {\rm R}),$ the so called vertical rotation set can be
defined as follows, see \cite{misiu}: 
\begin{equation}
\label{rotset}\rho _V(\widehat{f})=\bigcap_{{
\begin{array}{c}
i\geq 1 \\ 
\end{array}
}}\overline{\bigcup_{{
\begin{array}{c}
n\geq i \\ 
\end{array}
}}\left\{ \frac{p_2\circ \widehat{f}^n(\widehat{z})-p_2(\widehat{z})}n:
\widehat{z}\in S^1\times {\rm I}\negthinspace {\rm R}\right\} }
\end{equation}


This set is a closed interval (maybe a single point, but never empty) and it
was proved in \cite{eu1} and \cite{eu4} (and much earlier in \cite{doeff},
although the first author discovered this only recently) that all numbers in
its interior are realized by compact $f$-invariant subsets of ${\rm T^2,}$
which are periodic orbits in the rational case. From its definition, it is
easy to see that 
$$
\rho _V(\widehat{f}^m+(0,n))=m.\rho _V(\widehat{f})+n\text{ for any integers 
}n,m. 
$$

\item  Given $f\in DT({\rm T^2})$ and a lift $\widehat{f}\in DT(S^1\times 
{\rm I}\negthinspace {\rm R}),$ let $\mu $ be a $f$-invariant Borel
probability measure. We define the vertical rotation number of $\mu $ as
follows:%
$$
\rho _V(\mu )=\int_{{\rm T^2}}\phi (x,y)d\mu , 
$$
where the vertical displacement function $\phi :{\rm T^2}\rightarrow {\rm I}%
\negthinspace {\rm R}$ is given by $\phi (x,y)=p_2\circ \widehat{f}(\widehat{%
x},\widehat{y})-\widehat{y},$ for any $(\widehat{x},\widehat{y})\in
S^1\times {\rm I}\negthinspace {\rm R}$ such that $\pi ^{-1}(\widehat{x},
\widehat{y})\subset p^{-1}(x,y).$
\end{enumerate}

So, given $f\in DT({\rm T^2})$ and $\widehat{f}\in DT(S^1\times {\rm I}%
\negthinspace {\rm R}),$ as we said above, one wants to know, under which
conditions can $f$ be minimal? It is not difficult to see that in this case
the vertical rotation interval must be a single point, otherwise there would
be infinitely many periodic orbits. But more can be said.

\begin{theorem}
\label{irracional}: Given $f\in DT({\rm T^2})$ and a lift $\widehat{f}\in
DT(S^1\times {\rm I}\negthinspace {\rm R}),$ suppose that $f$ is minimal.
Then, $\rho _V(\widehat{f})=\{\alpha \}$ for some irrational number $\alpha .
$ \ 
\end{theorem}

So, if $f$ is a $C^r$ diffeomorphism, for some $r\geq 2,$ is there a natural
perturbation that destroys minimality?\ As the extreme points of $\rho _V( 
\widehat{f})$ vary continuously with $\widehat{f}\in DT(S^1\times {\rm I}%
\negthinspace {\rm R})$ (see \cite{doeff2}), a way to attack this problem is
by showing that irrational extremes are not stable under perturbations. This
was done in \cite{eu3} for twist mappings on the torus.

The main problem addressed in this paper is in a way, complementary to the
above. Suppose for instance that $\rho _V(\widehat{f})$ contains a single
rational number $p/q.$ What can we say about the dynamics of $f?$ And if $f$
preserves area and the center of gravity, that is Lebesgue measure has zero
vertical rotation number, what can we say about its vertical rotation
interval? When it is not reduced to zero, is zero always an interior point?
This is the so called Boyland's Conjecture.

Below we state our main results:

\begin{theorem}
\label{main1}: Given $f\in DT({\rm T^2})$ and a lift $\widehat{f}\in
DT(S^1\times {\rm I}\negthinspace {\rm R}),$ if $\rho _V(\widehat{f}%
)=\{p/q\},$ for some rational $p/q,$ then there exists a compact connected
set $K\subset S^1\times {\rm I}\negthinspace {\rm R,}$ invariant under $
\widehat{f}^q-(0,p),$ which separates the ends of the cylinder. So, all
points have uniformly bounded orbits under the action of $\widehat{f}%
^q-(0,p).$
\end{theorem}

Note that no area preservation hypothesis appear in our theorem. The
following corollary is almost immediate:

\begin{corollary}
\label{corol1}: Given $f\in DT({\rm T^2})$ and a lift $\widehat{f}\in
DT(S^1\times {\rm I}\negthinspace {\rm R}),$ suppose that $\rho _V(\widehat{f%
})=[a,p/q],$ for some rational $p/q$ and some real $a$ smaller than $p/q.$
Then there exists $M>0$ such that for all $\widehat{z}\in S^1\times {\rm I}%
\negthinspace {\rm R,}$ $p_2\circ \widehat{f}^n(\widehat{z})-p_2(\widehat{z}%
)-np/q<M,$ for all integers $n>0.$
\end{corollary}

The next result gives an explicit criteria which implies non-degenerate
vertical rotation sets and thus by a result analogous to the one in \cite
{llibre}, implies positive topological entropy (see for instance \cite{doeff}
and \cite{eu1}).

\begin{theorem}
\label{corol2}: Given $f\in DT({\rm T^2})$ and a lift $\widehat{f}\in
DT(S^1\times {\rm I}\negthinspace {\rm R}),$ there exists $M>0$ (which can
be explicitly computed) such that if for some points $\widehat{z}_1,\widehat{%
z}_2\in S^1\times {\rm I}\negthinspace {\rm R,}$ we have $p_2\circ \widehat{f%
}^{n_1}(\widehat{z}_1)-p_2(\widehat{z}_1)<-M$ and $p_2\circ \widehat{f}%
^{n_2}(\widehat{z}_2)-p_2(\widehat{z}_2)>M,$ for certain positive integers $%
n_1$ and $n_2,$ then $0$ is an interior point of $\rho _V(\widehat{f}).$
\end{theorem}

The next result gives a positive answer for Boyland's conjecture in this
setting:

\begin{corollary}
\label{main2}: Given an area-preserving $f\in DT({\rm T^2})$ and a lift $
\widehat{f}\in DT(S^1\times {\rm I}\negthinspace {\rm R})$ with zero
Lebesgue measure vertical rotation number, then either $\rho _V(\widehat{f})$
is reduced to $0,$ or $0$ is an interior point of $\rho _V(\widehat{f}).$
\end{corollary}

This paper is organized as follows. In the second section we present some
background results we use, with references and a few proofs and in the third
section we prove our main results. From now on we assume, without loss of
generality, that any $f\in DT({\rm T^2})$ we consider, is homotopic to a
Dehn twist $(x,y)\longrightarrow (x+k_{Dehn}y$ mod $1,y$ mod $1)$ with $%
k_{Dehn}>0.$

\section{Basic Tools}

\subsection{Brick Decompositions of the plane}

We define a brick decomposition of the plane as follows:

$$
{\rm I}\negthinspace {\rm R^2}=\stackrel{\infty }{\stackunder{i=0}{\cup }}%
D_i, 
$$
where each $D_i\in Brick\_Decomposition$ is the closure of a connected
simply connected open set, such that $\partial D_i$ is a polygonal simple
curve and $interior(D_i)\cap interior(D_j)=\emptyset ,$ for $i\neq j.$
Moreover, the decomposition is locally finite, that is, $\stackrel{\infty }{%
\stackunder{i=0}{\cup }}\partial D_i$ is a graph whose vertices have three
edges adjacent to them and the number of elements of the decomposition
contained in any compact subset of the plane is finite.

Given an orientation preserving homeomorphism of the plane $\widetilde{h},$
we say that the brick decomposition is free, if all its bricks are free,
that is, $\widetilde{h}(D_i)\cap D_i=\emptyset ,$ for all $i\in {\rm I}%
\negthinspace {\rm N.\ }$Given two bricks, $D$ and $E,$we say that there is
a chain connecting them, if there are bricks%
$$
D=D_0,\ D_1,D_2,...,D_{n-1},D_n=E 
$$
such that $\widetilde{h}(D_i)\cap D_{i+1}\neq \emptyset ,$ for $%
i=0,1,...,n-1.$ If $D=E,$ the chain is said to be closed.

In the following we will present a version of a theorem of J. Franks \cite
{franksanal} due to Le Roux and Guillou, see \cite{fleur}, page 39:

\begin{lemma}
\label{rouxgui}: The existence of a closed chain of free closed bricks
implies that there exists a simple closed curve $\gamma \subset {\rm I}%
\negthinspace {\rm R^2,}$ such that 
$$
index(\gamma ,\widetilde{h})=degree(\gamma ,\frac{\widetilde{h}(z)-z}{%
\left\| \widetilde{h}(z)-z\right\| })=1. 
$$
\end{lemma}

This result is a clever application of Brouwer's lemma on translation arcs.

\subsection{On the sets B$_S^{-}$ and B$_N^{+}$}

Here we present a theory developed in \cite{eufa1} and extend some
constructions to our new setting. For this, consider a homeomorphism $f\in
DT({\rm T^2}),$ a lift $\widehat{f}\in DT(S^1\times {\rm I}\negthinspace 
{\rm R})$ and a lift of $\widehat{f}$ to the plane, denoted $\widetilde{f}%
\in DT({\rm I\negthinspace R^2}).$ Given a real number $a$, let 
$$
H_a=S^1\times \{a\}, 
$$
$$
H_a^{-}=S^1\times ]-\infty ,a]\text{ and }H_a^{+}=S^1\times [a,+\infty [. 
$$

We will also denote the sets $H_0$, $H_0^{-}$ and $H_0^{+}$ simply by $H$, $%
H^{-}$ and $H^{+}$ respectively. If we consider the closed sets,

$$
B^{-}=\stackunder{n\leq 0}{\bigcap }_{}\widehat{f}^n({H}^{-}) 
$$
$$
and 
$$
$$
B^{+}=\stackunder{n\leq 0}{\bigcap }_{}\widehat{f}^n({H}^{+}), 
$$
we get that they are both closed and positively $\widehat{f}$-invariant. For
each of these sets, consider the following subsets: $B_S^{-}\subset B^{-}$
and $B_N^{+}\subset B^{+},$ each of which consisting of exactly all
unbounded connected components of respectively, $B^{-}$ and $B^{+}.$ The
sets $B_S^{-}$ and $B_N^{+}$ are always closed, but in some cases may be
empty. The next lemma tells us that under certain conditions, they really
exist.

\begin{lemma}
\label{existbgeral}:\ Suppose $0\in \rho _V(\widehat{f}).$ Then $B_N^{+}$
and $B_S^{-}$ are not empty.
\end{lemma}

{\it Proof:}

The proof of this result goes back to Le Calvez \cite{lecalvez2} and even
Birkhoff \cite{birk}.

First, suppose that $\stackunder{n\geq 0}{\cup }\widehat{f}^n(H)$ is
unbounded both from above and from below. In this case, considering the set $%
B_S^{-},$ the only thing we have to prove is that, for all $a\leq -1,$ there
exists a first positive integer $n=n(a),$ such that 
\begin{equation}
\label{onlymiss}\widehat{f}^{-n}(H_a)\cap H\neq \emptyset \text{ and }%
n(a)\rightarrow \infty \text{ as }a\rightarrow -\infty .
\end{equation}
Our assumption on $\stackunder{n\geq 0}{\cup }\widehat{f}^n(H)$ implies that 
$\widehat{f}^N(H)\cap H_a^{-}\neq \emptyset $ for some integer $N>0.$ If 
expression (\ref{onlymiss}) does not hold for $N,$ then $\widehat{f}%
^{-N}(H_a)\subset H^{+}\subset H_a^{+}+(0,1),$ which would imply that $%
0\notin \rho _V(\widehat{f}),$ a contradiction. So expression (\ref{onlymiss}%
) is true and the proof continues, for instance as in lemma 6 of \cite{eufa1}%
. A similar argument holds for $B_N^{+}$ (in this case $a\geq 1$).


If for some integer $M_0>0,$ $\widehat{f}^n(S^1\times \{0\})\subset
S^1\times [-M_0,+\infty [$ for all integers $n\geq 0,$ then clearly $
\widehat{f}^n(S^1\times [M_0,+\infty [)\subset S^1\times [0,+\infty [$ for
all integers $n\geq 0,$ so $B_N^{+}\supset S^1\times [M_0,+\infty [$ and
thus, it is not empty.

To prove that $B_S^{-}$ is also not empty, we have to work a little more.%

Let $O^{*}=\stackunder{n\geq 0}{\cup }\widehat{f}^n(S^1\times ]0,+\infty [)$
and let $O$ be the complement of the connected component of $(O^{*})^c$
which contains the lower end of the cylinder. We claim that $O^c$ is
connected and the same holds for $\partial $$O\stackrel{def.}{=}K.$ This
follows if we consider the $North-South$ compactification of the cylinder
and remember that it is a classical result, in the plane or sphere, that the
frontier of any connected component of the complement of a compact connected
subset is also connected. Clearly, $O^{*}\subset O$ (we just fill the
holes), $O$ contains the upper end of the cylinder and $\widehat{f}%
(O)\subset O.$

If $\stackunder{n\leq 0}{\cap }\widehat{f}^n(O^c)=\emptyset ,$ then $0\notin
\rho _V(\widehat{f}).$ So $\stackunder{n\leq 0}{\cap }\widehat{f}^n(O^c)\neq
\emptyset $ and as each connected component of this closed $\widehat{f}$%
-invariant set is bounded from above and unbounded, we get that for a
sufficiently large integer $j\geq 0,$ $\stackunder{n\leq 0}{\cap }\widehat{f}%
^n(O^c)-(0,j)\subset B_S^{-}\neq \emptyset .$

The remaining possibility can be treated in an analogous way. $\Box $

\vskip0.2truecm

\subsection{The $\omega $-limit sets of B$_S^{-}$ and B$_N^{+}$}

In this subsection we examine some properties of the set 
\begin{equation}
\label{defomelim}\omega (B_S^{-})\stackrel{def.}{=}\bigcap_{n=0}^\infty 
\overline{\bigcup_{i=n}^\infty \widehat{f}^i(B_S^{-})}. 
\end{equation}

Due to the fact that $\widehat{f}(B_S^{-})\subset B_S^{-}=\overline{B_S^{-}}%
, $ we get that 
\begin{equation}
\label{eelevadopi}\omega (B_S^{-})=\bigcap_{n=0}^\infty \widehat{f}%
^n(B_S^{-})=\bigcap_{n=-\infty }^\infty \widehat{f}^n(B_S^{-}). 
\end{equation}

\begin{lemma}
\label{propomegalim}: $\omega (B_S^{-})$ is a closed, $\widehat{f}$%
-invariant set, whose connected components are all unbounded.
\end{lemma}

{\it Proof:}

See lemma 7 of \cite{eufa1}. $\Box $

\vskip 0.2truecm

So from (\ref{eelevadopi}), $\omega $$(B_S^{-})\subset B_S^{-}$ and it is
still possible that $\omega (B_S^{-})=\emptyset .$ The next lemma tells us
that in this case, things are easier.

\begin{lemma}
\label{omelim}:\ Suppose $0\in \rho _V(\widehat{f}),$ which implies that $%
B_S^{-}$ is not empty. If $\omega (B_S^{-})=\emptyset ,$ then $\rho _V(
\widehat{f})\supset [-\epsilon ,0],$ for some $\epsilon >0.$
\end{lemma}

{\it Proof:}

See the proof of lemma 10 of \cite{eufa1} and the paragraph below it. $\Box $

\vskip 0.2truecm

Now, if we consider the set $B_S^{-}$ for $\widehat{f}^{-1},$ denoted $%
B_S^{-}(inv),$ we get the following:

\begin{lemma}
\label{omeliminv}: The sets $\omega (B_S^{-})$ and $\omega (B_S^{-}(inv))$
are equal.
\end{lemma}

{\it Proof: }

Let $\Gamma $ be a connected component of $\omega (B_S^{-}).$ From the
definition, $\widehat{f}^n(\Gamma )\subset H^{-}$ for all integers $n.$ So $%
\Gamma \subset B_S^{-}(inv)$ and moreover, for each positive integer $n,$ as 
$\widehat{f}^n(\Gamma )$ is contained in $H^{-},$ we get that $\Gamma
\subset \widehat{f}^{-n}(B_S^{-}(inv)),$ which means that $\Gamma \subset
\omega (B_S^{-}(inv)). $ Thus $\omega (B_S^{-})\subset \omega
(B_S^{-}(inv)). $ The other inclusion is proved in an analogous way. $\Box $

\vskip 0.2truecm

The following are important results on the structure of these sets.

\begin{lemma}
\label{ilim}: Any connected component $\widetilde{\Gamma }$ of $\pi ^{-1}($$%
\omega (B_S^{-}))$ is unbounded, not necessarily in the $\widetilde{y}$%
-direction.
\end{lemma}

{\it Proof:}

Let $d$ be the metric on $S^1\times $ ${\rm I}\negthinspace {\rm R}$ and let 
$\widetilde{d}$ be the lifted metric on the plane. Consider a point $
\widetilde{P}\in \widetilde{\Gamma }$ and let $P=\pi (\widetilde{P}).$ As $%
P\in \omega (B_S^{-}),$ there exists a connected component $\Gamma $ of $%
\omega (B_S^{-})$ that contains $P.$ Since by lemma \ref{propomegalim} $%
\Gamma $ is unbounded, for every sufficiently large integer $n$ there exists
a simple continuous arc $\gamma _n\subset S^1\times $ ${\rm I}%
\negthinspace 
{\rm R}$ such that:

\begin{itemize}
\item  $P$ is one endpoint of $\gamma _n;$

\item  $\gamma _n$ is contained in $S^1\times [-n,0]$ and it intersects $%
S^1\times \{-n\}$ only at its other endpoint;

\item  $\gamma _n$ is contained in a $(1/n,d)$-neighborhood of $\Gamma ;$
\end{itemize}

Now let $\widetilde{\gamma }_n$ be the connected component of $\pi
^{-1}(\gamma _n$$)$ that contains $\widetilde{P}.$ This arc $\widetilde{%
\gamma }_n$ is contained in a $(1/n,\widetilde{d})$-neighborhood of $\pi
^{-1}($$\Gamma )\subset \pi ^{-1}(\omega (B_S^{-}))$ because the covering
map is locally an isometry$.$

Now, embed the plane in the sphere $S^2=$ ${\rm I}\negthinspace {\rm %
R^2\sqcup \{\infty \}}$ equipped with a metric $D$ topologically equivalent
to the metric $\widetilde{d}$ on the plane. Then there exists a subsequence $
\widetilde{\gamma }_{n_i}\stackrel{i\rightarrow \infty }{\rightarrow }\Theta 
$ in the Hausdorff topology, for some compact connected set $\Theta \subset
S^2$. Clearly, both $\infty $ and $\widetilde{P}$ belong to $\Theta .$
Furthermore, since $\pi ^{-1}($$\omega (B_S^{-}))\cup \{\infty \}$ is a
closed set and 
$$
\stackunder{n\rightarrow \infty }{\lim }\left( \stackunder{\widetilde{z}\in 
\widetilde{\gamma }_n}{\sup }\text{ }\widetilde{d}(\widetilde{z},\pi
^{-1}(\omega (B_S^{-})))\right) =0, 
$$
we get that $\pi ^{-1}(\omega (B_S^{-}))\cup \{\infty \}$ contains $\Theta $
and the proof is over. 
%
$\Box $

\vskip 0.2truecm

\begin{lemma}
\label{complemen}: For any connected component $\widetilde{\Gamma }$ of $\pi
^{-1}($$\omega (B_S^{-})),$ $\widetilde{\Gamma }^c$ is connected.
\end{lemma}

{\it Proof: }

Take a connected component $\widetilde{\Gamma }$ of $\pi ^{-1}($$\omega
(B_S^{-})).$ First note that $\widetilde{\Gamma }^c$ has one connected
component, denoted $O^{+},$ which contains ${\rm I}\negthinspace {\rm R}%
\times ]0,+\infty [.$ So, if there is another one, denoted $O_1,$ it must be
contained in ${\rm I}\negthinspace {\rm R}\times ]-\infty ,0].$ In the
following we will prove that $\widetilde{f}^n(O_1)\subset {\rm I}%
\negthinspace {\rm R}\times ]-\infty ,0]$ for all integers $n.$

By contradiction, suppose that 
\begin{equation}
\label{refref1}\text{there is an integer }n_0\text{ such that }\widetilde{f}%
^{n_0}(O_1)\text{ is not contained in }{\rm I}\negthinspace 
{\rm R}\times ]-\infty ,0]. 
\end{equation}

There exists a number $m_0>0$ such that if $\widetilde{y}>m_0,$ then the
point $\widetilde{f}^{-n_0}(\widetilde{x},\widetilde{y})$ has positive $
\widetilde{y}$-coordinate, for all $\widetilde{x}\in {\rm I}\negthinspace 
{\rm R}$ (see (\ref{afbf}))${\rm .}$ So our hypothesis in (\ref{refref1})
implies that $\widetilde{f}^{-n_0}({\rm I}\negthinspace {\rm R}\times
]0,\infty [)\cap \partial O_1\neq \emptyset ,$ which means that $\widetilde{f%
}^{n_0}(\partial O_1)$ intersects ${\rm I}\negthinspace {\rm R}\times
]0,\infty [,$ a contradiction with the fact that 
$$
\widetilde{f}^{n_0}(\partial O_1)\subset \widetilde{f}^{n_0}(\widetilde{%
\Gamma })\subset \pi ^{-1}(\omega (B_S^{-}))\subset {\rm I}\negthinspace 
{\rm R}\times ]-\infty ,0]. 
$$
So (\ref{refref1}) does not hold. To conclude, let $\Gamma $ be the
connected component of $\omega (B_S^{-})$ that contains $\pi (\widetilde{%
\Gamma }),$ which as we know by lemma \ref{propomegalim} is unbounded. The
set $O_1\cup \widetilde{\Gamma }$ is connected as well as $\pi (O_1\cup 
\widetilde{\Gamma })\cup \Gamma =\pi (O_1)\cup \Gamma $ and the later is
contained in $\omega (B_S^{-})$ because $\widehat{f}^n(\pi (O_1))\subset
H^{-}$ for all integers $n.$ It follows that $\pi (O_1)\cup \Gamma =\Gamma
\subset \omega (B_S^{-})$ and therefore $O_1\cup \widetilde{\Gamma }$ is
contained in $\pi ^{-1}(\omega ($$B_S^{-})),$a contradiction with the choice
of $\widetilde{\Gamma }.$ $\Box $

\vskip 0.2truecm %
%

Clearly, similar results hold for $B_N^{+}.$

\section{Proofs}

\subsection{Proof of theorem 1}

Assume $f\in DT({\rm T^2})$ and its lift $\widehat{f}\in DT(S^1\times {\rm I}%
\negthinspace {\rm R})$ are such that $f$ is minimal and $\rho _V(\widehat{f}%
)$ is rational. Without loss of generality we can assume that $\rho _V( 
\widehat{f})=0,$ because if $f$ is minimal, the same happens for all its
iterates. This follows from the fact that, if for some integer $q>0,$ $f^q$
is not minimal, then it has a compact invariant minimal set $K\subset {\rm %
T^2},$ which by minimality, has empty interior. But then, 
$$
K\cup f(K)\cup ...\cup f^{q-1}(K) 
$$
is invariant under $f$ and as $K^c$ is open and dense, Baire 's property
also implies that $K\cup f(K)\cup ...\cup f^{q-1}(K)$ has empty interior, a
contradiction with the minimality of $f.$

As $\widehat{f}$ has no fixed points, lemma 2 of \cite{eu4} implies that
there exists a homotopically non trivial simple closed curve $\gamma $ in
the cylinder such that $\gamma \cap \widehat{f}(\gamma )=\emptyset .$
Without loss of generality, we can suppose that $\widehat{f}(\gamma )\subset
\gamma ^{-},$ the connected component of $\gamma ^c$ which is below $\gamma
. $ Let $k>0$ be an integer such that $\gamma -(0,k)\subset \gamma ^{-}.$ If
for some $n>0,$ $\widehat{f}^n(\gamma )\subset \left( \gamma -(0,k)\right)
^{-},$ then $0$ would not belong to $\rho _V(\widehat{f}).$ So, for all $%
n>0, $ there exists a point $\widehat{z}_n,$ above $\widehat{f}(\gamma )$
and below $\gamma ,$ such that 
$$
\{\widehat{z}_n,\widehat{f}(\widehat{z}_n),\widehat{f}^2(\widehat{z}_n),..., 
\widehat{f}^n(\widehat{z}_n)\}\text{ is above }\gamma -(0,k). 
$$
Taking a subsequence if necessary, we can assume that $\widehat{z}_{n_i}%
\stackrel{i\rightarrow \infty }{\rightarrow }\widehat{z}^{*},$ a point in
the closure of the region between $\widehat{f}(\gamma )$ and $\gamma .$
Clearly, the positive orbit of $\widehat{z}^{*}$ is bounded in the cylinder
and so its $\omega $-limit set $\omega (\widehat{z}^{*})$ is a compact $
\widehat{f}$-invariant subset of the cylinder. Moreover, as any integer
vertical translate of $\omega (\widehat{z}^{*})$ is also $\widehat{f}$%
-invariant, if we pick a minimal $\widehat{f}$-invariant compact set $K$
contained in $\omega (\widehat{z}^{*}),$ clearly, by minimality it satisfies 
$K\cap K+(0,n)=\emptyset $ for all $n\neq 0.$

As $f$ is minimal, when $K$ is projected to the torus is must be the whole
torus, a contradiction. $\Box $

\subsection{Proof of theorem 2}

Given $f\in DT({\rm T^2})$ and a lift $\widehat{f}\in DT(S^1\times {\rm I}%
\negthinspace {\rm R}),$ without any loss of generality we can assume that $%
\rho _V(\widehat{f})=0.$

Lemma \ref{existbgeral} implies that $B_N^{+}\neq \emptyset $ and $%
B_S^{-}\neq \emptyset ,$ and lemma \ref{omelim} implies that the same holds
for their $\omega $-limits, $\omega (B_N^{+})\neq \emptyset $ and $\omega
(B_S^{-})\neq \emptyset .$

In the following we will present two technical results. For each $\widehat{x}%
\in S^1,$ consider the following functions, which as the next lemma shows,
are well defined at all $\widehat{x}\in S^1:$%
$$
\begin{array}{c}
\mu ( 
\widehat{x})=\max \{\widehat{y}\in {\rm I}\negthinspace {\rm R:}\text{ }( 
\widehat{x},\widehat{y})\in \omega (B_S^{-})\} \\ \nu (\widehat{x})=\min \{ 
\widehat{y}\in {\rm I}\negthinspace {\rm R:}\text{ }(\widehat{x},\widehat{y}%
)\in \omega (B_N^{+})\} 
\end{array}
$$

\begin{lemma}
\label{unifomega}: There exists a constant $M_f>0$ such that 
$$
\stackunder{\widehat{x},\widehat{y}\in S^1}{\sup }\left| \mu (\widehat{x}%
)-\mu (\widehat{y})\right| \leq M_f\text{ and }\stackunder{\widehat{x},
\widehat{y}\in S^1}{\sup }\left| \nu (\widehat{x})-\nu (\widehat{y})\right|
\leq M_f.\text{ } 
$$
\end{lemma}

{\it Proof:}

The proof is analogous for both cases, so let us only consider the function $%
\mu .$ As $\omega (B_S^{-})$ is closed and bounded from above, choose some $
\widehat{x}_0\in S^1$ such that $\{\widehat{x}_0\}\times ]-\infty ,0]\cap
\omega (B_S^{-})\neq \emptyset $ and for some $\widehat{y}_0\leq 0,$ $(
\widehat{x}_0,\widehat{y}_0)$ belongs to $\omega (B_S^{-})$ and has maximal $
\widehat{y}$-coordinate. Then $\mu (\widehat{x}_0)=\widehat{y}_0$ is well
defined.

Note that as $f$ is homotopic to a Dehn twist, for all $(\widetilde{x}, 
\widetilde{y})\in ${\rm I}\negthinspace 
{\rm R$^2$} there are constants $A_f>0$ and $B_f>0$ such that 
\begin{equation}
\label{afbf}\left| p_2\circ \widetilde{f}(\widetilde{x},\widetilde{y})- 
\widetilde{y}\right| <A_f\text{ and }\left| p_1\circ \widetilde{f}( 
\widetilde{x},\widetilde{y})-\widetilde{x}-k_{Dehn}\widetilde{y}\right|
<B_f. \text{ } 
\end{equation}

So for any compact set $G\subset {\rm I}\negthinspace {\rm R^2}$ with

$$
\left| p_2(G)\right| \stackrel{def.}{=}\max (p_2(G))-\min (p_2(G))\geq V_f%
\stackrel{def.}{=}\frac{(3+2B_f)}{k_{Dehn}} 
$$
$$
\text{and} 
$$
$$
\left| p_1(G)\right| \stackrel{def.}{=}\max (p_1(G))-\min (p_1(G))<1, 
$$

we have:

$$
\left| p_1(\widetilde{f}(G))\right| >2 \text{ and } p_2\mid _{\widetilde{f}%
(G)}>\min (p_2(G))-A_f. 
$$

Consider the intersection $\pi ^{-1}(\omega (B_S^{-}))\cap {\rm I}%
\negthinspace {\rm R}\times [\mu (\widehat{x}_0)-V_f,\mu (\widehat{x}_0)].$
If all vertical segments $Seg_{\widetilde{x}}=\{\widetilde{x}\}\times [\mu (
\widehat{x}_0)-V_f,\mu (\widehat{x}_0)]$ intersect $\pi ^{-1}(\omega
(B_S^{-})),$ then for all $\widehat{x}\in S^1,$ $\mu (\widehat{x}_0)-V_f\leq
\mu (\widehat{x})\leq 0$ and the proof is over. So, suppose that there
exists a real number $\widetilde{x}^{*}$ such that $Seg_{\widetilde{x}^{*}}$
do not intersect $\pi ^{-1}(\omega (B_S^{-})).$ This implies that for any
integer $n,$ $Seg_{\widetilde{x}^{*}}+(n,0)$ do not intersect $\pi
^{-1}(\omega (B_S^{-})).$ Let $\theta $ be the connected component of $%
\omega (B_S^{-})$ containing $(\widehat{x}_0,\widehat{y}_0)$ and let $\Theta 
$ be a component of $\pi ^{-1}(\theta ).$ The set $\Theta $ is also a
connected component of $\pi ^{-1}(\omega (B_S^{-})),$ so by lemma \ref{ilim}
it is unbounded. It is now clear that $\Theta $ intersects the two
horizontal boundaries of $[\widetilde{x}^{*}+n_\Theta ,\widetilde{x}%
^{*}+n_\Theta +1]\times [\mu (\widehat{x}_0)-V_f,\mu (\widehat{x}_0)]$ for
some integer $n_\Theta ,$ because it can not meet the open half plane $\{
\widetilde{y}>\mu (\widehat{x}_0)\}.$

Thus, $\left| p_1(\widetilde{f}(\Theta ))\right| >2$ and $p_2\mid _{
\widetilde{f}(\Theta )}>\mu (\widehat{x}_0)-V_f-A_f.$ As $\omega (B_S^{-})$
is invariant, $\pi \left( \widetilde{f}(\Theta )\right) $$\subset \omega
(B_S^{-})$ and so for any $\widehat{x}\in S^1,$ $\mu (\widehat{x}%
_0)-V_f-A_f<\mu (\widehat{x})\leq 0.$

The above argument implies that if we choose $M_f=V_f+A_f,$ then we are
done. $\Box $

\vskip 0.2truecm

Now, let us define the number 
\begin{equation}
\label{defmdehnw}M_{Dehn}=\frac{2+B_f}{k_{Dehn}}>0. 
\end{equation}
A simple computation shows that for all $(\widetilde{x},\widetilde{y})\in 
{\rm I}\negthinspace {\rm R^2}$ with $\widetilde{y}>M_{Dehn},$ we have

$$
p_1\circ \widetilde{f}(\widetilde{x},\widetilde{y})>\widetilde{x}+2\text{
and }p_1\circ \widetilde{f}(\widetilde{x},-\widetilde{y})<\widetilde{x}-2. 
$$

The construction performed below is analogous for both $\omega (B_N^{+})$
and $\omega (B_S^{-}).$ The details will be presented for $\omega (B_S^{-}).$
First, note that for every $\widehat{x}\in S^1,$ $\mu (\widehat{x})+\left( -%
\stackunder{\widehat{z}\in S^1}{\max }\mu (\widehat{z})+M_f\right)
+M_{Dehn}\geq M_{Dehn}.$ This means that if we define the following positive
integer number $n_{trans}\stackrel{def.}{=}$ $\left\lfloor -\stackunder{ 
\widehat{z}\in S^1}{\max }\mu (\widehat{z})+M_f+M_{Dehn}\right\rfloor +1$ ($%
\left\lfloor a\right\rfloor $ is the integer part of $a$), then the set 
\begin{equation}
\label{ometrans}\omega (B_S^{-})_{trans}\stackrel{def.}{=}\omega
(B_S^{-})+\left( 0,n_{trans}\right) 
\end{equation}
has, for every $\widehat{x}\in S^1,$ a point of the form $(\widehat{x}, 
\widehat{y}),$ with $\widehat{y}>M_{Dehn}.$ In other words, the function $%
\mu _{trans}$ associated with $\omega (B_S^{-})_{trans}$ satisfies $\mu
_{trans}(\widehat{x})\stackrel{def.}{=}\mu (\widehat{x})+n_{trans}>M_{Dehn},$
for all $\widehat{x}\in S^1.$

Now, for a fixed $\widetilde{x}\in {\rm I}\negthinspace {\rm R,}$ consider
the semi-line $\{\widetilde{x}\}\times [M_{Dehn},+\infty [.$ When we
intersect it with 
$$
\widetilde{\omega (B_S^{-})}_{trans}\stackrel{def.}{=}\pi ^{-1}\left( \omega
(B_S^{-})_{trans}\right) 
$$
we get that $\{\widetilde{x}\}\times ]\mu _{trans}(\pi (\widetilde{x}%
)),+\infty [\cap \widetilde{\omega (B_S^{-})}_{trans}=\emptyset $ (note that 
$\widetilde{\omega (B_S^{-})}_{trans}$ is also a $\widetilde{f}$-invariant
set).

Let $v=\{\widetilde{x}\}\times ]\mu _{trans}(\pi (\widetilde{x})),+\infty [$
and let $\Theta $ be the connected component of $\widetilde{\omega (B_S^{-})}%
_{trans}$ that contains $(\widetilde{x},\mu _{trans}(\pi (\widetilde{x}))).$

\begin{lemma}
\label{tetav}: The following holds: $\Theta \cup v$ is a closed connected
set, $\left( \Theta \cup v\right) ^c$ has two open connected components, one
of which is positively invariant and $\widetilde{f}^n(v)\cap v=\emptyset $
for all integers $n\neq 0.$
\end{lemma}

{\it Proof:}

The fact that $\Theta \cup v$ is closed and connected is obvious. As $\Theta 
$ is a connected component of $\widetilde{\omega (B_S^{-})}_{trans},$ it is
unbounded and limited from above in the $\widetilde{y}$-direction. ${\rm \ }$

By the Jordan separation theorem, we get that $\left( \Theta \cup v\right)
^c $ has at least two connected components, $O_L$ and $O_R,$ defined as
follows: For any point $\widetilde{P}\in v,$ there exists $\delta $$>0$ such
that $B_\delta (\widetilde{P})\cap \Theta =\emptyset .$ Moreover, $B_\delta
( \widetilde{P})\backslash v$ has exactly 2 connected components, one to the
left of $v,$ contained in $O_L$ and the other one to the right of $v,$
contained in $O_R.$ So their closures, $\overline{O_L}$ and $\overline{O_R}$
both contain $v.$ Now, suppose $\left( \Theta \cup v\right) ^c$ has another
connected component, denoted $O^{*}.$ Clearly $\partial $$O^{*}$ do not
intersect $v$ because all points sufficiently close to a point in $v$ and,
not in $v,$ are contained in $O_L\cup O_R.$ So, $\partial $$O^{*}\subset
\Theta $ and $O^{*}$ is then a connected component of $\Theta ^c$ bounded
from above in the $\widetilde{y}$-direction. And this contradicts lemma \ref
{complemen}. So, $\left( \Theta \cup v\right) ^c=O_L\cup O_R.$

Note that $\widetilde{f}(v)\cap v=\widetilde{f}(v)\cap \Theta =\widetilde{f}%
^{-1}(v)\cap \Theta =\emptyset .$ The paragraph after definition (\ref
{defmdehnw}) implies that $\widetilde{f}(v)\subset O_R.$ In the following we
will show that $\widetilde{f}(O_R)\subset O_R.$

There are 2 possibilities:

\begin{enumerate}
\item  {$\tilde f(\Theta )\neq \Theta $} $\Rightarrow $ {$\tilde f(\Theta
)\cap \Theta =\emptyset ,$} because $\Theta $ is a connected component of an
invariant set;

\item  {$\tilde f(\Theta )=\Theta ;$}
\end{enumerate}

Assume first that $\tilde f(\Theta )\cap \Theta =\emptyset $. Then 
$$
\widetilde{f}(\Theta \cup v)\cap (\Theta \cup v)=\emptyset . 
$$
Since $\widetilde{f}(v)\subset O_R$ and $\widetilde{f}(\Theta \cup v)$ is
connected, we get that $\widetilde{f}(\Theta \cup v)\subset O_R,$ so $%
O_L\cup \Theta \cup v$ is contained either in $\widetilde{f}(O_L)$ or $
\widetilde{f}(O_R).$ It can not be contained in $\widetilde{f}(O_R)$ because
a point of the form $(-a,a)$ for a sufficiently large $a>0$ is contained in $%
O_L$ and $\widetilde{f}^{-1}(-a,a)$ is also contained in $O_L,$ see (\ref
{afbf}). Thus, $O_L\cup \Theta \cup v\subset \widetilde{f}(O_L),$ which
implies that, $\widetilde{f}(O_R)\subset O_R.$

Now suppose $\tilde f(\Theta )=\Theta .$ This implies that $O_L\cup v\subset
(\widetilde{f}(v\cup \Theta ))^c$ because $\widetilde{f}(v)\subset O_R$ and $
\widetilde{f}(\Theta )=\Theta .$ So, by connectedness, $O_L\cup v$ is
contained either in $\widetilde{f}(O_R)$ or in $\widetilde{f}(O_L).$ As in
the case $\tilde f(\Theta )\cap \Theta =\emptyset ,$ one actually gets $%
O_L\cup v\subset \widetilde{f}(O_L)$ so%
$$
\widetilde{f}(O_R)\subset (\widetilde{f}(O_L))^c\subset (O_L\cup
v)^c=O_R\cup \Theta  
$$
and since $\widetilde{f}(\Theta )=\Theta ,$ we finally get that $\widetilde{f%
}(O_R)\subset O_R.$

In order to finish the proof, note that, as $\widetilde{f}(v)\cap
v=\emptyset ,$ for any $n\geq 2,$ $\widetilde{f}^n(v)\subset \widetilde{f}%
(O_R),$ which do not intersect $v.$ So $\widetilde{f}^n(v)\cap v=\emptyset .$
This finishes the proof of our lemma. $\Box $

\vskip 0.2truecm

{\bf Remarks}:

\begin{itemize}
\item  as $\mu _{trans}(\pi (\widetilde{x}))<M_f+M_{Dehn}+2$ for all $
\widetilde{x}\in ${\rm I}\negthinspace {\rm R}, we get that $\widetilde{f}%
^n(\{\widetilde{x}\}\times [M_f+M_{Dehn}+2,+\infty [)\cap \{\widetilde{x}%
\}\times [M_f+M_{Dehn}+2,+\infty [=\emptyset $ for all integers $n>0.$

\item  an analogous argument applied to $\omega (B_N^{+})$ implies that for
any $\widetilde{x}\in ${\rm I}\negthinspace {\rm R, }if $w=\{\widetilde{x}%
\}\times ]-\infty ,\nu (\pi (\widetilde{x}))-\left\lfloor \stackunder{
\widehat{z}\in S^1}{\inf }\nu (\widehat{z})+M_f+M_{Dehn}\right\rfloor -1[,$
then $\widetilde{f}^n(w)\cap w=\emptyset $ for all integers $n>0.$ So as in
the above remark, $\nu _{trans}(\pi (\widetilde{x}))>-2-M_f-M_{Dehn}$ for
all $\widetilde{x}\in ${\rm I}\negthinspace {\rm R}, which implies that $
\widetilde{f}^n(\{\widetilde{x}\}\times ]-\infty ,-M_f-M_{Dehn}-2[)\cap \{
\widetilde{x}\}\times ]-\infty ,-M_f-M_{Dehn}-2[=\emptyset $ for all
integers $n>0.$
\end{itemize}

\vskip 0.2truecm

Summarizing, there exists a real number $M^{\prime }>0$ such that for all $
\widetilde{x}\in ${\rm I}\negthinspace {\rm R}, $\widetilde{f}^n(\{ 
\widetilde{x}\}\times [M^{\prime },+\infty [)\cap \{\widetilde{x}\}\times
[M^{\prime },+\infty [=\emptyset $ and $\widetilde{f}^n(\{\widetilde{x}%
\}\times ]-\infty ,-M^{\prime }])\cap \{\widetilde{x}\}\times ]-\infty
,-M^{\prime }]=\emptyset $ for all integers $n>0$ and 
\begin{equation}
\label{defmprima}M^{\prime }\stackrel{def.}{=}M_f+M_{Dehn}+2=\frac{5+3B_f}{%
k_{Dehn}}+A_f+2 
\end{equation}

Now let us suppose by contradiction that there exists a point $\widehat{z}$
in the cylinder and an integer $n_0>0$ such that%
$$
\left| p_2(\widehat{f}^{n_0}(\widehat{z}))-p_2(\widehat{z})\right|
>2M^{\prime }+8. 
$$
Without loss of generality, we can assume that $p_2(\widehat{z})<-M^{\prime
}-3$ and $p_2(\widehat{f}^{n_0}(\widehat{z}))>M^{\prime }+3.$

Let us also consider the fixed point free mapping of the plane 
$$
\widetilde{g}(\bullet )=\widetilde{f}^{n_0}(\bullet )-(0,1). 
$$
To see that it is actually fixed point free, note that if $\widetilde{g}$
has a fixed point, then $1/n_0\in \rho _V(\widehat{f}),$ a contradiction.
Now, note that for all $\widetilde{x}\in ${\rm I}\negthinspace 
{\rm R, }$\widetilde{g}(\{\widetilde{x}\}\times [M^{\prime }+2,+\infty
[)\cap \{\widetilde{x}\}\times [M^{\prime }+2,+\infty [=\emptyset $ and $
\widetilde{g}(\{\widetilde{x}\}\times ]-\infty ,-M^{\prime }-2])\cap \{ 
\widetilde{x}\}\times ]-\infty ,-M^{\prime }-2]=\emptyset .$ Moreover, using
the fact that $\widetilde{g}$ is also the lift of a torus homeomorphism
homotopic to a Dehn twist and a compacity argument, one can prove that there
exists an integer $N>0,$ such that for all integers $n,$ the sets 
\begin{equation}
\label{deffenes} 
\begin{array}{c}
F_n^{-}=[n/N,(n+1)/N]\times ]-\infty ,-M^{\prime }-2] \\ 
\text{and} \\ F_n^{+}=[n/N,(n+1)/N]\times [M^{\prime }+2,\infty [ 
\end{array}
\end{equation}
are free under $\widetilde{g},$ that is, $\widetilde{g}(F_n^{+or-})\cap
F_n^{+or-}=\emptyset ,$ for all integers $n.$ Moreover, the fact that $%
k_{Dehn}>0$ (see the end of section 1) implies that there exists an integer $%
K_{crit}>0,$ such that for all integers $n$

$$
\begin{array}{c}
\widetilde{g}(F_n^{+})\cap F_m^{+}\neq \emptyset ,\text{ for all }m\geq
n+K_{crit} \\ \text{and} \\ \widetilde{g}(F_n^{-})\cap F_m^{-}\neq \emptyset
,\text{ for all }m\leq n-K_{crit}. 
\end{array}
$$

These will be important bricks in a special brick decomposition of the plane
in $\widetilde{g}$-free sets we will construct, which will be invariant
under integer horizontal translations $(\widetilde{x},\widetilde{y}%
)\rightarrow (\widetilde{x}+1,\widetilde{y}).$

Clearly, such a construction is possible, because as $\widetilde{g}( 
\widetilde{x}+1,\widetilde{y})=\widetilde{g}(\widetilde{x},\widetilde{y}%
)+(1,0),$ we just have to decompose $S^1\times [-M^{\prime }-2,M^{\prime
}+2] $ into a union of bricks with sufficiently small diameter, so that
their pre-images under $\pi $ are $\widetilde{g}$-free.

To conclude our proof, we will show that this brick decomposition has a
closed brick chain, a contradiction with the fact that $\widetilde{g}$ is
fixed point free, see lemma \ref{rouxgui}. This idea was already used in the
proof of theorem 4 of \cite{eu4}.

Consider a point $\widetilde{z}\in \pi ^{-1}(\widehat{z})$ and a brick $%
F_{i_0}^{-}$ that contains $\widetilde{z}.$ From our choices, 
$$
\widetilde{g}(F_{i_0}^{-})\cap F_{i_1}^{+}\neq \emptyset ,\text{ for some
integer }i_1. 
$$

As $\rho _V(\widehat{f})=\{0\},$ let us choose a point $\widehat{w}\in
S^1\times ]M^{\prime }+2,+\infty [${\rm \ }such that%
$$
p_2(\widehat{g}^n(\widehat{w}))\stackrel{n\rightarrow \infty }{\rightarrow }%
-\infty , 
$$
where $\widehat{g}(\bullet )\stackrel{def.}{=}\widehat{f}^{n_0}(\bullet
)-(0,1)$ (as $\rho _V(\widehat{g})=\{-1\},$ all points in $S^1\times {\rm I}%
\negthinspace 
{\rm R\ }$satisfy the above condition). So, we can choose a point $
\widetilde{w}\in F_{i_2}^{+},$ for some integer $i_2,$ such that:

\begin{itemize}
\item  $i_2>i_1+K_{crit},$ so $\widetilde{g}(F_{i_1}^{+})\cap
F_{i_2}^{+}\neq \emptyset ;$

\item  $\widetilde{g}^{n_2}(\widetilde{w})\in F_{i_3}^{-},$ for some
integers $n_2>0$ and $i_3>i_0+K_{crit};$
\end{itemize}

As $\widetilde{g}(F_{i_3}^{-})\cap F_{i_0}^{-}\neq \emptyset ,$ we get there
exists a closed brick chain starting at $F_{i_0}^{-}.$ As we said, this is a
contradiction because $\widetilde{g}$ is fixed point free. Thus $\widehat{f}%
^n(S^1\times \{0\})\subset S^1\times [-8-2M^{\prime },2M^{\prime }+8]$ for
all integers $n>0.$ In order to conclude the proof, let $K$ be the only
connected component of the frontier of 
$$
\stackunder{n\geq 0}{\cap }\widehat{f}^n(closure(\stackunder{m\geq 0}{\cup } 
\widehat{f}^m(S^1\times ]0,+\infty [))) 
$$
which does not bound a disc. Then $K$ is a compact connected set that
separates the ends of the cylinder, $\widehat{f}(K+(0,l))=K+(0,l),$ for all
integers $l$ and $\left| p_2(K)\right| \leq 4M^{\prime }+20.$ $\Box $

\subsection{Proof of corollary 1}

Without loss of generality, by considering $\widehat{f}^q-(0,p),$ we can
suppose that $\rho _V(\widehat{f})=[a,0],$ for some $a<0.$ As in the proof
of theorem 2, lemma \ref{existbgeral} implies that $B_N^{+}\neq \emptyset ,$ 
$B_S^{-}\neq \emptyset $ and $B_N^{+}(inv)\neq \emptyset ,$ $%
B_S^{-}(inv)\neq \emptyset .$ If for instance $\omega (B_S^{-})=\emptyset ,$
then lemma \ref{omeliminv} implies that $\omega (B_S^{-}(inv))=\emptyset $
and so lemma \ref{omelim} implies that there exists $\epsilon >0$ such that $%
\rho _V(\widehat{f}^{-1})\supset [-\epsilon ,0],$ which gives $\rho _V( 
\widehat{f})\supset [0,\epsilon ],$ a contradiction$.$ So, we can assume
that $\omega (B_N^{+})\neq \emptyset $ and $\omega (B_S^{-})\neq \emptyset .$
If we suppose that for every $M>0,$ there exists a point $\widehat{z}\in
S^1\times {\rm I}\negthinspace {\rm R}$ and an integer $n>0$ such that 
$$
p_2(\widehat{f}^n(\widehat{z}))-p_2(\widehat{z})>M, 
$$
then following exactly the same ideas used in theorem 2, we arrive at a
contradiction which proves the corollary. $\Box $

\subsection{Proof of theorem 3}

As in theorem 2, let us fix a $\widetilde{f}\in DT({\rm I\negthinspace R^2}%
), $ which is a lift of $\widehat{f}.$ First, we will show that if 
$$
M\geq M_0\stackrel{def.}{=}(20+2B_f)/k_{Dehn}+10\text{ (see (\ref{afbf})),} 
$$
then $\widehat{f}$ has a fixed point. In case $\widehat{f}$ is fixed point
free, lemma 2 of \cite{eu4} tells us that there exists a homotopically
non-trivial simple closed curve $\gamma \subset S^1\times {\rm I}%
\negthinspace {\rm R}$ such that $\widehat{f}(\gamma )\cap \gamma =\emptyset 
$ and $\gamma \subset S^1\times [-m_D,m_D],$ where $m_D>0$ is the smallest
real number that satisfies 
\begin{equation}
\label{defmmm} 
\begin{array}{c}
\widetilde{f}(\{\widetilde{x}\}\times [m_D,+\infty [)\subset [\widetilde{x}%
+10,+\infty [\times {\rm I}\negthinspace {\rm R}\text{ } \\ and \\ 
\text{ }\widetilde{f}(\{\widetilde{x}\}\times [-\infty ,-m_D])\subset
]-\infty ,\widetilde{x}-10]\times {\rm I}\negthinspace {\rm R,} 
\end{array}
\end{equation}
for all $\widetilde{x}\in {\rm I\negthinspace R.}$ A simple computation
shows that if we take $m_D$ equal $(10+B_f)/k_{Dehn},$ then (\ref{defmmm})
is satisfied.

So, as $M\geq 2m_D+10,$ the theorem hypotheses imply that $\widehat{f}$ has
a fixed point. Thus $0\in \rho _V(\widehat{f})$ and lemma \ref{existbgeral}
implies that $B_N^{+}\neq \emptyset ,$ $B_S^{-}\neq \emptyset $ and the same
holds for the inverse of $\widehat{f},$ namely, $B_S^{-}(inv)\neq \emptyset $
and $B_N^{+}(inv)\neq \emptyset .$ If $\omega (B_N^{+})=\emptyset ,$ then
lemma \ref{omelim} implies that there exists $\delta >0$ such that $\rho _V( 
\widehat{f})\supset [0,\delta ].$ Also, from lemma \ref{omeliminv} we get
that $\omega (B_N^{+}(inv))=\emptyset $ and so again by lemma \ref{omelim},
there exists $\epsilon >0$ such that $\rho _V(\widehat{f}^{-1})\supset
[0,\epsilon ],$ which gives $\rho _V(\widehat{f})\supset [-\epsilon ,\delta
] $ and the theorem is proved. So, again we can suppose that $\omega
(B_S^{-})\neq \emptyset $ and $\omega (B_N^{+})\neq \emptyset .$

If $\rho _V(\widehat{f})=[a,0]$ for some $a\leq 0,$ then if 
$$
M\geq M_1\stackrel{def.}{=}2M^{\prime }+8=\frac{10+6B_f}{k_{Dehn}}+2A_f+12, 
$$
by the same argument used to prove theorem 2, we arrive at a contradiction.
The same happens in the other possibility, that is, if $\rho _V(\widehat{f}%
)=[0,b],$ for some $b>0.$

So, it is enough to choose 
$$
M=\max \{M_0,M_1\}\leq \frac{20+6B_f}{k_{Dehn}}+2A_f+12\text{ to finish the
proof}.\text{ }\Box 
$$

\subsection{Proof of Corollary 2}

Let us start by showing that there are two possibilities:

1) $\stackunder{n\geq 0}{\cup }\widehat{f}^n(H)$ is bounded and this means
that $\rho _V(\widehat{f})=\{0\};$

2) $\stackunder{n\geq 0}{\cup }\widehat{f}^n(H)$ is unbounded from above and
from below;

In order to understand that the above are the only possible cases, suppose
for instance that $\stackunder{n\geq 0}{\cup }\widehat{f}^n(H)$ is unbounded
and contained in $H_a^{+}$ for some real number $a<0.$

As in lemma \ref{existbgeral}, let $O^{*}=\stackunder{n\geq 0}{\cup } 
\widehat{f}^n(S^1\times ]0,+\infty [)$ and let $O$ be the complement of the
connected component of $(O^{*})^c$ which contains the lower end of the
cylinder. As in that lemma, $\partial $$O\stackrel{def.}{=}K$ is a compact
connected set that separates the ends of the cylinder. Clearly, $%
O^{*}\subset O$ (we just fill the holes), $H_1^{+}\subset O\subset H_a^{+},$ 
$O$ is an open set homeomorphic to the cylinder and $\widehat{f}(O)\subset
O. $

Let us state a simple result, but before we present a definition:

\begin{description}
\item[Definition]  : If $\gamma $ is a homotopically non trivial simple
closed curve in $S^1\times {\rm I}\negthinspace {\rm R,}$ then $\gamma ^c%
\stackrel{def.}{=}\gamma ^{-o}\cup \gamma ^{+o},$ where $\gamma ^{-o(+o)}$
is the open connected component of $\gamma ^c$ which contains the lower
(upper) end of the cylinder. We define $\gamma ^{-}\stackrel{def.}{=}%
closure(\gamma ^{-o})=\gamma ^{-o}\cup \gamma $ and the same for $\gamma
^{+}.$
\end{description}

\begin{proposition}
\label{exata}: Given an area-preserving $f\in DT({\rm T^2})$ and a lift $
\widehat{f}\in DT(S^1\times {\rm I}\negthinspace {\rm R})$ with zero
Lebesgue measure vertical rotation number, for any $b\in {\rm I}%
\negthinspace {\rm R}$ the following equality holds (in this case $\widehat{f%
}$ is said to be exact):%
$$
Leb(H_b^{+}\cap (\widehat{f}(H_b))^{-})=Leb(H_b^{-}\cap (\widehat{f}%
(H_b))^{+}), 
$$
where for any measurable set $D,$ $Leb(D)\stackrel{def.}{=}Lebesgue$ $measure
$ of $D.$
\end{proposition}

{\it Proof:}

If we remember (\ref{afbf}), we get that there exists an integer $N>0$ such
that, for any given $b\in {\rm I}\negthinspace {\rm R,\ }\widehat{f}%
(H_b)\cap (H_{b+N}\cup H_{b-N})=\emptyset .$ So, consider the finite annulus 
$\Omega \stackrel{def.}{=}$$S^1\times [b,b+N].$ As it is a finite union of
fundamental domains of the torus, we get that 
\begin{equation}
\label{fabiohoje}\int_\Omega \left[ p_2\circ \widehat{f}(\widehat{x}, 
\widehat{y})-\widehat{y}\right] d\widehat{x}d\widehat{y}=0\text{ (this
follows from }\rho _V(Leb)=0\text{).}\ 
\end{equation}

Note that we can write

$$
\Omega =\left( \widehat{f}(\Omega )\cap \Omega \right) \cup \left(
H_b^{+}\cap (\widehat{f}(H_b))^{-o}\right) \cup \left( H_b^{-}\cap (\widehat{%
f}(H_b))^{+o}+(0,N)\right) 
$$

$$
and 
$$
$$
\widehat{f}(\Omega )=\left( \widehat{f}(\Omega )\cap \Omega \right) \cup
\left( H_b^{+o}\cap (\widehat{f}(H_b))^{-}+(0,N)\right) \cup \left(
H_b^{-o}\cap (\widehat{f}(H_b))^{+}\right) , 
$$
where the unions are disjoint$.$ Expression (\ref{fabiohoje}) together with
the preservation of area imply that the $\widehat{y}$-coordinate of the
geometric center of $\Omega $ and of $\widehat{f}(\Omega )$ are equal. So,
let us compute them (for a measurable set $\Pi $ in the cylinder, we denote
the $\widehat{y}$-coordinate of its geometric center by $\widehat{y}%
_{G.C.(\Pi )})$:%
$$
\widehat{y}_{G.C.(\Omega )}=\left[ 
\begin{array}{c}
\widehat{y}_{G.C.(\widehat{f}(\Omega )\cap \Omega )}.Leb(\widehat{f}(\Omega
)\cap \Omega )+ \\ +
\widehat{y}_{G.C.(H_b^{+}\cap (\widehat{f}(H_b))^{-})}.Leb(H_b^{+}\cap (
\widehat{f}(H_b))^{-})+ \\ +\left( \widehat{y}_{G.C.(H_b^{-}\cap (\widehat{f}%
(H_b))+)}+N\right) .Leb(H_b^{-}\cap (\widehat{f}(H_b))^{+})
\end{array}
\right] /Leb(\Omega ) 
$$

$$
\widehat{y}_{G.C.(\widehat{f}(\Omega ))}=\left[ 
\begin{array}{c}
\widehat{y}_{G.C.(\widehat{f}(\Omega )\cap \Omega )}.Leb(\widehat{f}(\Omega
)\cap \Omega )+ \\ +\left( 
\widehat{y}_{G.C.(H_b^{+}\cap (\widehat{f}(H_b))^{-})}+N\right)
.Leb(H_b^{+}\cap (\widehat{f}(H_b))^{-})+ \\ +\widehat{y}_{G.C.(H_b^{-}\cap
( \widehat{f}(H_b))+)}.Leb(H_b^{-}\cap (\widehat{f}(H_b))^{+}) 
\end{array}
\right] /Leb(\widehat{f}(\Omega )) 
$$
As $Leb(\widehat{f}(\Omega ))=Leb(\Omega )$ and $\widehat{y}_{G.C.(\widehat{f%
}(\Omega ))}=\widehat{y}_{G.C.(\Omega )},$ we get that 
$$
N.Leb(H_b^{+}\cap (\widehat{f}(H_b))^{-})=N.Leb(H_b^{-}\cap (\widehat{f}%
(H_b))^{+}), 
$$
which proves the proposition (note that we used the fact that $Leb(H_b)=0$). 
$\Box $

\vskip 0.2truecm

Now let us choose $c\in {\rm I}\negthinspace {\rm R}$ such that $\{K\cup 
\widehat{f}(K)\}\subset interior(H_c^{-}\cap (\widehat{f}(H_c))^{-}).$ From
the preservation of Lebesgue measure and the above proposition, we get that 
$$
Leb(O\cap H_c^{-})=Leb(\widehat{f}(O)\cap (\widehat{f}(H_c))^{-})=Leb(
\widehat{f}(O)\cap H_c^{-}). 
$$
The choice of $c,$ together with the fact that $\widehat{f}(O)\subset O,$
implies that $closure(O)=closure(\widehat{f}(O))=\widehat{f}(closure(O)).$
So $\partial (closure(O))$ separates the ends of the cylinder and is $
\widehat{f}$-invariant. But this means that all orbits are uniformly
bounded, a contradiction with our hypothesis that $\stackunder{n\geq 0}{\cup 
}\widehat{f}^n(H)$ is unbounded. So, either 1) or 2) from the beginning of
the proof of the corollary can happen.

And in possibility 2) we can apply theorem 3 to conclude the proof. $\Box $


\begin{thebibliography}{99}
\bibitem{eu1}  Addas-Zanata S. (2002): On the existence of a new type of
periodic and quasi-periodic orbits for twist maps of torus. {\it Nonlinearity%
} {\bf 15}, 1399-1416.

\bibitem{eu3}  Addas-Zanata S. (2005): On properties of the vertical
rotation interval for twist mappings. {\it Erg. Th. \& Dyn. Sys.} {\bf 25},
641-660

\bibitem{eu4}  Addas-Zanata S. (2005): Some extensions of the
Poincar\'e-Birkhoff theorem to the cylinder and a remark on mappings of the
torus homotopic to Dehn twists. {\it Nonlinearity} {\bf 18}, 2243-2260.

\bibitem{eufa1}  Addas-Zanata S. and Tal F. A. (2009): Homeomorphisms of the
annulus with a transitive lift. {\it Math. Zeit.} {\bf 267}, 971-980.

\bibitem{birk}  Birkhoff G. (1935): Nouvelles Recherches sur les syst\`emes
dynamiques. Collec. Math. Pap. {\bf 2,} 530-661.

\bibitem{doeff}  Doeff E. (1997): Rotation measures for homeomorphisms of
the torus homotopic to a Dehn twist. {\it Erg. Th. \& Dyn. Sys.} {\bf 17},
1-17.

\bibitem{doeff2}  Doeff E. and Misiurewicz M. (1997): Shear rotation
numbers. {\it Nonlinearity} {\bf 10}, 1755-1762.

\bibitem{franksanal}  Franks J. (1988): Generalizations of the
Poincar\'e-Birkhoff theorem. {\it Ann. of Math. }{\bf 128, }139-151

\bibitem{fleur}  Le Roux F.(2004): Hom\'eomorphismes de surfaces:
th\'eor\`emes de la fleur de Leau-Fatou et de la vari\'et\'e stable. {\it %
Ast\'erisque} {\bf 292}, 210 pp.

\bibitem{lecalvez2}  Le Calvez P. (1987) Propri\'et\'es dynamiques des
r\'egions d'instabilit\'e. {\it Ann. Sci. \'Ecole Norm. Sup.} {\bf 20},
443-464.

\bibitem{llibre}  Llibre J. and Mackay R. (1991): Rotation vectors and
entropy for homeomorphisms of the torus isotopic to the identity. {\it \
Erg. Th. \& Dyn. Sys. }{\bf 11, }115-128.

\bibitem{misiu}  Misiurewicz M. and Ziemian K. (1989): Rotation sets for
maps of tori. {\it J. London Math. Soc.} {\bf 40} 490-506.

\end{thebibliography}
\end{document}